\documentclass[a4paper,12pt]{article}
\usepackage{a4wide}
\usepackage{amsmath}
\usepackage{amssymb}
\usepackage{amsthm}
\usepackage{latexsym}
\usepackage{graphicx}
\usepackage[english]{babel}
\usepackage{makeidx}
              
\newtheorem{dfn} [subsection]{Definition}
\newtheorem{obs} [subsection]{Remark}
\newtheorem{exm} [subsection]{Example}

\newtheorem{prop}[subsection]{Proposition}

\newtheorem{teor}[subsection]{Theorem}
\newtheorem{lema}[subsection]{Lemma}
\newtheorem{cor} [subsection]{Corollary}

\newcommand{\de}{\mathbf{d}}
\newcommand{\me}{\mathbf{m}}
\newcommand{\en}{\mathbf{n}}

\newcommand{\di}{$\de$-fixed ideal }
\newcommand{\dis}{$\de$-fixed ideals }

\begin{document}

\selectlanguage{english}
\frenchspacing

\large
\begin{center}
\textbf{Regularity for certain classes of monomial ideals.}

Mircea Cimpoea\c s
\end{center}

\normalsize

\footnotetext[1]{This paper was supported by the CEEX Program of the Romanian
Ministry of Education and Research, Contract CEX05-D11-11/2005 and by 
 the Higher Education Commission of Pakistan.}

\begin{abstract}
We introduce a new class of monomial ideals, called strong Borel type ideals, and we compute the Mumford-Castelnouvo regularity for principal strong Borel type ideals. Also, we describe the $\de$-fixed ideals generated by powers of variables and we compute their regularity.

\vspace{5 pt} \noindent \textbf{Keywords:} p-Borel ideals, Borel type ideals, Mumford-Castelnuovo regularity.

\vspace{5 pt} \noindent \textbf{2000 Mathematics Subject
Classification:}Primary: 13P10, Secondary: 13E10.
\end{abstract}

\begin{center}
\textbf{Introduction.}
\end{center}

Let $K$ be an infinite field, and let $S=K[x_1,...,x_n],n\geq 2$ the polynomial ring over $K$.
Bayer and Stillman \cite{bs} note that a Borel fixed ideal $I$ satisfies the following property $(I:x_j^\infty)=(I:(x_1,\ldots,x_j)^\infty)$ for all $j=1,\ldots,n$. Herzog, Popescu and Vladoiu say that
a monomial ideal is of Borel type if it fulfill the previous condition. We mention that this concept appears also in \cite[Definition 1.3]{cs} as the so called weakly stable ideal. In fact, Herzog, Popescu and Vladoiu notice that a monomial ideal $I$ is of Borel type, if and only if for any monomial $u\in I$ and for any $1\leq j<i \leq n$, there exists an integer $t>0$ such that $x_j^{t}u/x_i^{\nu_i(u)}\in I$, where $\nu_i(u)>0$ is the exponent of $x_i$ in $u$. (see \cite[Proposition 1.2]{hpv}). This property suggest us to define the so called ideals of strong Borel type (Definition $1.1$), or simply, (SBT)-ideals. In the first section, we give the explicit form of a principal (SBT)-ideal (Lemma $1.4$) and we compute its regularity (Theorem $1.6$).

Let $\de: 1=d_{0}|d_{1}|\cdots|d_{s}$ be a strictly increasing sequence of positive integers. We say that
$\de$ is a $\de$-sequence. In \cite{mir} it was proved that for any $a\in \mathbb N$ there exists an unique sequence of positive integers $a_{0},a_{1},\ldots,a_{s}$ such that: $a= \sum_{t=0}^{s}a_{t}d_{t}$ and $0 \leq a_{t} < \frac{d_{t+1}}{d_{t}}$, for any $0 \leq t<s$.
The decomposition $a=\sum_{t=0}^{s}a_{t}d_{t}$ is called the $\de$-decomposition of $a$. In particular, if $d_{t}=p^{t}$ we get the $p$-adic decomposition of $a$. Let $a,b\in\mathbb N$ and consider the decompositions $a=\sum_{t=0}^{s}a_{t}d_{t}$ and $b=\sum_{t=0}^{s}b_{t}d_{t}$. We say that $a\leq_{\de}b$ if $a_{t}\leq b_{t}$ for any $0\leq t\leq s$. We say that a monomial ideal $I\subset S$ is $\de$-fixed, if for any monomial $u\in I$ and for any indices $1\leq j<i \leq n$, if $t\leq_{\de} \nu_{i}(u)$ then $u \cdot x_{j}^{t}/x_{i}^{t}\in I$ (see \cite[Definition 1.4]{mir}).

In \cite{mir}, it was proved a formula for the regularity of a principal \di, i.e the smallest \di which contains a given
monomial $u\in S$. This formula generalize the Pardue's formula for the regularity of a principal $p$-Borel ideal,
proved in \cite{ah} and \cite{hp}, and later in \cite{hpv}. In the section $2$, we describe the $\de$-fixed ideals generated by powers of variables (Proposition $2.2$) and we give a formula for their regularity (Corollary $2.8$).

The author owes a special thanks to Assistant Professor Alin Stefan for valuable discussions on section $2$ of this paper. My thanks goes also to the School of Mathematical Sciences, GC University, Lahore, Pakistan for supporting 
and facilitating this research.
\section{Monomial ideals of strong Borel type.}

Let $K$ be an infinite field, and let $S=K[x_1,...,x_n],n\geq 2$ the polynomial ring over $K$.

\begin{dfn}
We say that a monomial ideal $I\subset S$ is of strong Borel type (SBT) if for any monomial $u\in I$ and for any $1\leq j<i \leq n$, there exists an integer $0\leq t\leq \nu_i(u)$ such that $x_j^{t}u/x_i^{\nu_i(u)}\in I$.
\end{dfn}

\begin{obs}
Obviously, an ideal of strong Borel type is also an ideal of Borel type, but the converse is not true. Take for instance
$I=(x_1^3,x_2^2)\subset K[x_1,x_2]$.

The sum of two ideals of (SBT) is still an ideal of (SBT). The same is true for
an intersection or a product of two ideals of (SBT).
\end{obs}

\begin{dfn}
Let $\mathcal A\subset S$ be a set of monomials. We say that $I$ is the (SBT)-ideal generated by $\mathcal A$, if $I$
is the smallest, with respect to inclusion, ideal of (SBT) containing $\mathcal A$. We write $I=SBT(\mathcal A)$.

In particular, if $\mathcal A=\{u\}$, where $u\in S$ is a monomial, we say that $I$ is the principal (SBT)-ideal
generated by $u$, and we write $I=SBT(u)$.
\end{dfn}

\begin{lema}
Let $1\leq i_1<i_2<\cdots <i_r\leq n$ be some integers, $\alpha_1,\ldots,\alpha_r$ some positive integers and
$u=x_{i_1}^{\alpha_1}x_{i_2}^{\alpha_2}\cdots x_{i_r}^{\alpha_r} \in S$. Then, the principal (SBT)-ideal generated by $u$, is:
\[ I = SBT(u) = \prod_{q=1}^{r}(\me_q^{[\alpha_q]}),\;where\;\me_q=\{x_1,\ldots,x_{i_q}\}\;and\; 
\me_q^{[\alpha_q]}=\{x_1^{\alpha_q},\ldots,x_{i_q}^{\alpha_q}\}. \]
\end{lema}

\begin{proof}
Denote $I'=\prod_{q=1}^{r}(\me_q^{[\alpha_q]})$. If $v$ is a minimal monomial generator of $I'$, then \linebreak
$v=x_{j_1}^{\alpha_1}x_{j_2}^{\alpha_2}\cdots x_{j_r}^{\alpha_r}$, for some $1\leq j_q\leq i_q$, where $1\leq q\leq r$. Since
\[ v = \frac{x_{j_r}^{\alpha_r}}{x_{i_r}^{\alpha_r}}\cdots \frac{x_{j_2}^{\alpha_2}}{x_{i_2}^{\alpha_2}}\cdot
\frac{x_{j_1}^{\alpha_1}}{x_{i_1}^{\alpha_1}} u, \]
and $I$ is of (SBT) it follows that $v\in I$ and thus $I' \subseteq I$. For the converse, simply notice that $I'$ is itself a (SBT)-ideal.
\end{proof}

\begin{obs}
For any monomial ideal $I\subset S$, we denote
$m(I)=max\{m(u):\;u\in G(I)\}$, where $G(I)$ is the set of the minimal generators of $I$ and $m(u)=max\{i:\;x_i|u\}$. Also, if $M$ is a graded $S$-module of finite length, we denote $s(M)=max\{t:\;M_{t}\neq 0\}$.

Let $I\subset S$ be a Borel type ideal. In \cite{hpv}, it is defined a chains of ideals $I=I_0 \subset I_1 \subset \cdots \subset I_r=S$ as follows. We let $I_0=I$. Suppose $I_{\ell}$ is already defined. If $I_{\ell}=S$ then the chain ends. Otherwise, we let $n_{\ell}=m(I_{\ell})$ and set
$I_{\ell+1}=(I_{\ell}:x_{n_{\ell}}^{\infty})$. Notice that $r\leq n$, since $n_{\ell}>n_{\ell+1}$ for all $0\leq \ell<r$.
The chain $I = I_{0} \subset I_{1}\subset \cdots \subset I_r = S$ is called the \emph{sequential chain} of $I$.
\cite[Corollary 2.5]{hpv} states that 
\[(1)\; I_{\ell+1}/I_{\ell} \cong (J_{\ell}^{sat}/J_{\ell})[x_{n_{\ell}+1},\ldots,x_n],\]
for all $0\leq \ell<r$, where $J_{\ell}\subset S_{\ell} = K[x_1,\ldots,x_{n_{\ell}}]$ is the ideal generated by $G(I_{\ell})$. Also,\linebreak \cite[Corollary 2.5]{hpv} gives a formula for the regularity of $I$, more precisely,
\[(2)\; reg(I) = max\{s(J_{0}^{sat}/J_0),s(J_{1}^{sat}/J_1),\cdots,s(J_{r-1}^{sat}/J_{r-1})\} + 1.\]
\end{obs}

Our next goal is to give a formula for the regularity of a principal (SBT)-ideal. In order to do it, we will use the
previous remark. 

Let $1\leq i_1<i_2<\cdots <i_r\leq n$ be some integers, $\alpha_1,\ldots,\alpha_r$ some positive integers and
$u=x_{i_1}^{\alpha_1}x_{i_2}^{\alpha_2}\cdots x_{i_r}^{\alpha_r} \in S$. For each $1\leq q\leq r$, $1\leq f\leq q$ with $\alpha_f\leq \alpha_q$ and $1\leq j\leq i_q$, we define the numbers:
\[ \chi_{qj}^{(f)}:=\begin{cases}
                      \alpha_j+\alpha_q - 1, & if \; j<q\;and\;\alpha_j\geq \alpha_f \\
                      \alpha_f - 1,& otherwise                
              \end{cases}, \chi_q^{(f)}:=\sum_{j=1}^{i_q}\chi_{qj}^{(f)} \; and\; \chi_q=\max_{f}\chi_q^{(f)}. \]

\begin{teor}
With the above notations, we have $reg(SBT(u))=\max\limits_{q=1}^{r}\chi_q + 1$.
\end{teor}

\begin{proof}
Firstly, we describe the sequential chain of $I$. Since $I_r:= I = \prod_{q=1}^{r}(\me_q^{[\alpha_q]})$, it follows that $I_{r-1} := (I_r:x_{i_r}^{\infty}) = \prod_{q=1}^{r-1}(\me_q^{[\alpha_q]})$. Analogously, we get $I_{q}:=(I_{q+1}:x_{i_{q+1}}^{\infty}) = \prod_{e=1}^{q}(\me_e^{[\alpha_e]})$, for all $0\leq q < r$. Therefore, the sequential chain of $I$ is, \[ I = I_r \subset I_{r-1} \subset \cdots \subset I_1 \subset I_0 = S. \]
Let $J_q$ be the ideal of $S_q = K[x_1,\ldots,x_{i_q}]$ generated by $G(I_q)$, for $1\leq q\leq r$. Denoting
$s_q = s(J_q^{sat}/J_q)$, (2) from Remark $1.5$ implies $reg(I)=max\{s_q:\; 1\leq q\leq r\}$, so, in order to compute the regularity of $I$,
we must determine the numbers $s_q$. We claim that $s_q = \chi_q$. 

First of all, note that $J_q = I_q\cap S_q$ and $J_q^{sat} = I_{q-1}\cap S_q$.
Let $1\leq f\leq q$ with $\alpha_f\leq \alpha_q$ and $w = x_1^{\chi_{q1}^{(f)}}\cdots x_{i_q}^{\chi_{q,i_q}^{(f)}}$. 
Since $\chi_{qe}^{(f)}\geq \alpha_e$ for any $1\leq e \leq q-1$ we get $x_1^{\chi_{q1}^{(f)}}\cdots x_{q-1}^{\chi_{q,q-1}^{(f)}}\in J_{q}^{sat} = \prod_{e=1}^{q-1}(\me_{e}^{[\alpha_e]})S_q$, therefore $w \in J_{q}^{sat}$. On the other hand, one can easily see that $w\notin J_q$, so $w$ is a nonzero element
in $J_{q}^{sat}/J_q$ with $deg(w)=\chi_q$, thus $s_q\geq \chi_q$. 

In order to prove the converse inequality, we consider a monomial $u\in J_{q}^{sat}$ with $deg(u)\geq \chi_q + 1$ and we show that $u\in J_q$. Assume by contradiction that $u\notin J_q$. Since $u\in J_q^{sat}$, it follows that
$u = x_{j_1}^{\alpha_1}\cdots x_{j_{q-1}}^{\alpha_{q-1}}\cdot x_1^{\beta_{1}}\cdots x_{i_q}^{\beta_{i_q}}$, where
$1\leq j_e \leq i_e$ for $1\leq e\leq q-1$ and $\beta_1 +\cdots + \beta_{i_q} \geq \chi_q - \sum_{e=1}^{q-1}\alpha_e$.
Let $A = \{1,\ldots,i_q\}\setminus \{j_1,\ldots,j_{q-1}\}$. Since $u\notin J_q$ and
$x_{j_1}^{\alpha_1}\cdots x_{j_{q-1}}^{\alpha_{q-1}}\in J_q^{sat}$ it follows $\beta_j\leq \alpha_q - 1$
for all $j\in A$.

Write $\{1,\ldots,q-1\} = \cup_{i=1}^{m}E_i$, where $E_i=\{e_{i1},\ldots,e_{ik_i}\}$, such that
$j_{e_{ik}}=j_{e_i}$ for all $1\leq k \leq k_i$ and $E_i\cap E_{i'} = \emptyset$ whenever $i\neq i'$. With these notations,
\[ u = x_{j_{e_1}}^{\alpha_{e_{11}} + \cdots + \alpha_{e_{1k_1}} + \beta_{j_{e_1}}}\cdots 
x_{j_{e_m}}^{\alpha_{e_{m1}} + \cdots + \alpha_{e_{mk_m}} + \beta_{j_{e_m}}}\cdot \prod_{j\in A}x_{j}^{\beta_j}.\]

Let $1\leq f\leq q$ such that $\alpha_f\leq \alpha_q$, $\beta_{j}<\alpha_f$ for all $j\in A$ and $\alpha_f$ is the largest integer between all the $\alpha_{f'}$, with $f'$ satisfying the above conditions. Suppose that there exist some
$1\leq i \leq m$ and $1\leq k\leq k_i$ such that $\alpha_{e_{ik}}<\alpha_q$. It follows that $\beta_{j_{e_i}}\leq \alpha_f - \alpha_{e_{ik}}-1$, otherwise $u\in J_q$. One can immediately conclude that 
$\sum_{e=1}^{q-1}\alpha_e + \sum_{j=1}^{i_q}\beta_j \leq \chi_{q}^{(f)}$.
\end{proof}

\begin{exm}{\em
 Let $u=x_2^{6}x_3^{7}\in S=K[x_1,x_2,x_3]$. From Lemma $1.4$ it follows that $I = SBT(u) =
	      (x_1^6,x_2^6)(x_1^7,x_2^7,x_3^7)$. With the notations of $1.5$ and $1.6$, we have $J_1=(x_1^6,x_2^6)\subset 
	      K[x_1,x_2]$ and $J_2 = I$. Also, $J_1^{sat}=K[x_1,x_2]$ and $J_2^{sat} = (x_1^6,x_2^6) \subset S$.
	      Obviously, $\chi_1=\chi_1^{(1)} = 2\cdot 5 = 10$, i.e. $s(J_1^{sat}/J_1) = s(K[x_1,x_2]/(x_1^6,x_2^6)) = 10$. 
	      We have $\chi_{2}^{(1)} = (6+7-1) + 2\cdot 5 = 23$ and $\chi_{2}^{(2)} = 3\cdot 6 = 18$, therefore
	      $\chi_2=23$ and thus $reg(I)=max\{10,23\} + 1 = 24$.
}\end{exm}

In the end of this section, we mention the following result, which generalize a result of Eisenbud-Reeves-Totaro (see \cite[Proposition 12]{ert}).

\begin{prop}\cite[Corollary 8]{mir2}
If $I$ is a Borel type ideal, then \[ reg(I)=min\{e:\;e\geq deg(I),\; I_{\geq e}\;is\; stable \},\]
where $deg(I)$ is the maximal degree of a minimal monomial generator of $I$.
\end{prop}

In particular, this holds for (SBT)-ideals, and thus we get the following corollary.

\begin{cor}
With the notations of Theorem $1.5$, if $I=SBT(u)$ and $e\geq \max\limits_{q=1}^{r}\chi_q + 1$ 
then $I_{\geq e}$ is stable.
\end{cor}

\begin{obs}
Note also that the regularity of a (SBT)-ideal, $I\subset S$, is upper bounded by $n(deg(I)-1)+1$, (see \cite[Theorem 2.2]{saf}). In fact, $deg(I)$ is the maximum degree of a minimal generator of $I$ as a (SBT)-ideal! 
\end{obs}

\section{\dis generated by powers of variables.}

Firstly, let fix some notations. Let $u_{1},\ldots,u_{m}\in S$ be some monomials. We say that $I$ is the \di generated by $u_{1},\ldots,u_{m}$, if $I$ is the smallest \di, w.r.t inclusion, which contain $u_{1},\ldots,u_{m}$, and we write $I=<u_{1},\ldots,u_{m}>_{\de}$. In particular, if $m=1$, we say that $I$ is the principal \di generated by $u=u_{1}$ and we write $I=<u>_{\de}$. 

In the case when $I$ is a principal \di, \cite[Theorem 3.1]{mir} gives a formula for the Castelnuovo-Mumford regularity of $I$. Using similar tehniques as in \cite{mir}, we will compute the regularity for \dis generated by powers of variables. We recall some results proved in \cite{mir} which are useful. Let $\alpha$ be a positive integer and let $I = <x_{n}^{\alpha}>_{\de} \subset S=K[x_{1},\ldots,x_{n}]$. Suppose $\alpha=\sum_{t=0}^{s}\alpha_{t}d_{t}$ with $\alpha_{s}\neq 0$. Then:
\begin{itemize}
	\item $I=\prod_{t=0}^{s} (\mathbf{m}^{[d_{t}]})^{\alpha_{t}}$, where $\mathbf{m}=\{x_{1},\ldots,x_{n}\}$ and
	      $\mathbf{m}^{[d]}=\{x_{1}^{d},\ldots,x_{n}^{d}\}$ \cite[1.6]{mir}.
	\item $Soc(S/I) = (J+I)/I$ with $J=\sum_{t=0}^{s}(x_{1}\cdots x_{n})^{d_{t}-1}(\mathbf{m}^{[d_{t}]})^{\alpha_{t}-1}
	      \prod_{j>t} (\mathbf{m}^{[d_{j}]})^{\alpha_{j}}$ \cite[2.1]{mir}.
	\item $reg(I) = max\{e:\; ((J+I)/I)_{e}\neq 0 \} = \alpha_{s}d_{s}+ (n-1)(d_{s}-1)$ (see \cite[3.1]{mir}).
	\item If $e\geq reg(I)$ then $I_{\geq e}$ is stable (see \cite[3.6]{mir} or apply Proposition $1.8$, since any $d$-fixed ideal is of Borel type, see \cite[1.11]{mir}).
\end{itemize}

\newpage
\begin{lema}
If $1\leq j\leq j'\leq n$ and $\alpha\geq \beta$ are positive integers, then $<x_{j}^{\alpha}>\subset <x_{j'}^{\beta}>$.
\end{lema}

\begin{proof}
Indeed, using \cite[1.7]{mir} it is enough to notice that $<x_{j}^{\alpha}> \subset <x_{j'}^{\alpha}>$, since $x_{j}^{\alpha} \in <x_{j'}^{\alpha}>$.
\end{proof}

Our next goal is to give the set of the minimal generators of a \di generated by some powers of variables. Using the previous lemma, we had reduced to the next case:

\begin{prop}
Let $n\geq 2$ and let $1\leq i_{1} < i_{2} < \cdots < i_{r}=n$ be some integers. Let $\alpha_{1}< \alpha_{2}< \cdots < \alpha_{r}$ be some positive integers. Then
\[ I=<x_{i_{1}}^{\alpha_{1}},x_{i_{2}}^{\alpha_{2}},\ldots,x_{i_{r}}^{\alpha_{r}}>_{\de} = \sum_{q=1}^{r} I^{(q)},
with \; I^{(q)} =
   \sum_{ \footnotesize \begin{array}{c} \gamma_{1},\ldots,\gamma_{q}\leq_{\de}\alpha_{q},\\ \;\gamma_{1}+\cdots+\gamma_{i}<\alpha_{i},\;for\;i<q\\ 
   \gamma_{1}+\cdots+\gamma_{i}<_{d}\alpha_{q},\;for\;i<q\\
   \gamma_{1}+\cdots+\gamma_{q}=\alpha_{q} \end{array} \normalsize }   
   \prod_{e=1}^{q}\prod_{t=0}^{s} (\mathbf{n}_{e}^{[d_{t}]})^{\gamma_{et}},\]
   where $\mathbf{n}_{e} = \{x_{i_{e-1}+1},\ldots,x_{i_{e}}\}$, $\mathbf{n}_{e}^{[d_{t}]} = \{x_{i_{e-1}+1}^{d_{t}}, \ldots,x_{i_{e}}^{d_{t}}\}$, $i_{0}=0$ and $\gamma_{e}=\sum_{t=0}^{s}\gamma_{et}d_t$.
\end{prop}

\begin{proof}
Let $\mathbf{m}_{q} = \{x_{1},\ldots,x_{i_{q}}\}$ for $1\leq q\leq r$. Obviously, $\en_{q}=\me_{q}\setminus \me_{q-1}$ for $q>1$ and $\me_{1}=\en_{1}$. Using the simple fact that $I$ is the sum of principal \dis generated by the $\de$-generators of $I$ together with \cite[Proposition 1.6]{mir} we get:
\[ I 
= \sum_{q=1}^{r} \prod_{t=0}^{s} (\me_{q}^{[d_{t}]})^{\alpha_{qt}},\;where\;\alpha_{q}=\sum_{t=0}^{s}\alpha_{qt}d_{t}  \]
Denote $S_{q}=K[x_{1},\ldots,x_{i_{q}}]$ for $1\leq q\leq r$. In order to obtain the required formula, we use induction on $r\geq 1$, the case $r=1$ being obvious. Let $r>1$ and
assume that the assertion is true for $r-1$, i.e 
\[ I'=<x_{i_{1}}^{\alpha_{1}},\ldots,x_{i_{r-1}}^{\alpha_{r-1}}>_{\de} = \sum_{q=1}^{r-1}  
   \sum_{\footnotesize \begin{array}{c} \gamma_{1},\ldots,\gamma_{q}\leq_{\de}\alpha_{q},\\ \;\gamma_{1}+\cdots+\gamma_{i}<\alpha_{i},\;for\;i<q\\
   \gamma_{1}+\cdots+\gamma_{i}<_{d}\alpha_{q},\;for\;i<q\\
   \gamma_{1}+\cdots+\gamma_{q}=\alpha_{q} \end{array} \normalsize} \prod_{e=1}^{q}\prod_{t=0}^{s} (\mathbf{n}_{e}^{[d_{t}]})^{\gamma_{et}}  \subset S_{r-1}. \]
Obviously, $I=I'S + <x_{n}^{\alpha_{r}}>_{\de} = I'S + \prod_{t=0}^{s}(\me_{r}^{[d_{t}]})^{\alpha_{rt}}$. Also, $I'S$ 
and $I'$ have the same set of minimal generators and none of the minimal generators of $I'S$ is in $I^{(r)}$.
But, a minimal generator of $<x_{n}^{\alpha_{r}}>_{\de}$ is of the form $w=\prod_{t=0}^{s}\prod_{j=1}^{n}x_{j}^{\lambda_{tj}d_{t}}$ with $0\leq \lambda_{tj}$ and $\sum_{j=1}^{n}\lambda_{tj}=\alpha_{rt}$. Suppose $w\notin I'S$. In order to complete the proof, we will 
show that $w\in I^{(r)}$. Let $v_{q}=\prod_{t=0}^{s}\prod_{j=i_{q-1}+1}^{i_{q}}x_{j}^{\lambda_{tj}d_{t}}$ and
let $w_{q}=\prod_{e=1}^{q}v_{e}$. Obvious, $w = v_{1}\cdots v_{r} = w_{r}$. Since $w\notin I'$ it follows that
$w_{q}\notin I^{(q)}$ for any $1\leq q\leq r-1$. But $w_{q}\notin I^{(q)}$ implies 
$(*)\;\sum_{t=0}^{s}\sum_{j=1}^{i_{q}} \lambda_{tj}d_{t} < \alpha_{q}$, otherwise $w_{q} \in <x_{i_{q}}^{\alpha_{q}}S_{q}>_{\de}S_{r-1} \subset I'$ and thus $w\in I'$, a contradiction.
We choose $\gamma_{e} = \sum_{t=0}^{s}\sum_{j=i_{e-1}+1}^{i_{e}}\lambda_{tj}d_{t}$ for $1\leq e\leq r$.
For $1\leq q<r$, $(*)$ implies $\gamma_{1}+\cdots +\gamma_{q}<\alpha_{q}$. On the other hand, it is obvious that 
$\gamma_{1}+\cdots + \gamma_{e}\leq_{d}\alpha_{r}$ for any $1\leq e\leq r$ and $\gamma_{1}+\cdots+\gamma_{r}=\alpha_{r}$. Thus $w\in I^{(r)}$ as required.
\end{proof}

\begin{exm}
Let $\de:1|2|4|12$ and let $I=<x_{2}^{7},x_{3}^{10},x_{5}^{17}>_{\de} \subset K[x_{1},\ldots,x_{5}]$. We have $7=1\cdot 1 + 1\cdot 2 + 1\cdot 4$, $10=1\cdot 2 + 2\cdot 4$, $17=1\cdot 1 + 1\cdot 4 + 1\cdot 12$. 
We have 

\[ I^{(1)}= <x_2^{7}>_{\de} = (x_{1},x_{2})(x_{1}^{2},x_{2}^{2})(x_{1}^{4},x_{2}^{4}).\]

In order to compute $I^{(2)}$, we need to find all the pairs $(\gamma_1,\gamma_2)$ such that $\gamma_1<7$,
$\gamma_1<_{\de}10$ and $\gamma_2=10-\gamma_1$. We have $4$ pairs, namely $(0,10)$, $(2,8)$, $(4,6)$ and $(6,4)$, thus

\[ I^{(2)}=(x_{1}^{2},x_{2}^{2})(x_{1}^{4},x_{2}^{4})x_{3}^{4} + (x_{1}^{4},x_{2}^{4})x_{3}^{6} + 
	      (x_{1}^{2},x_{2}^{2})x_{3}^{8} + (x_{3}^{10}).\]

In order to compute	$I^{(3)}$, we need to find all $(\gamma_1,\gamma_2,\gamma_3)$ such that $\gamma_1<7$,
$\gamma_1+\gamma_2<10$, $\gamma_1<_{\de}17$, $\gamma_1+\gamma_2<_{\de}17$ and $\gamma_3 = 17-\gamma_1+\gamma_2$.
If $\gamma_1=0$ then, the pair $(\gamma_2,\gamma_3)$ is one of the following:$(0,17)$,$(1,16)$,$(4,13)$ or $(5,12)$.
If $\gamma_1=1$ then, the pair $(\gamma_2,\gamma_3)$ is one of the following:$(0,16)$ of $(4,12)$.
If $\gamma_1=4$ then, the pair $(\gamma_2,\gamma_3)$ is one of the following:$(0,13)$ of $(1,12)$. If $\gamma_1=5$ then, the pair $(\gamma_2,\gamma_3)$ is $(0,12)$. Thus 
\[I^{(3)} = (x_{1},x_{2})(x_{1}^{4},x_{2}^{4})(x_{4}^{12},x_{5}^{12}) +
	      (x_{1}^{4},x_{2}^{4})x_{3}(x_{4}^{12},x_{5}^{12}) + (x_{1}^{4},x_{2}^{4})(x_{4},x_{5})(x_{4}^{12},x_{5}^{12})+ \]\[ + (x_{1},x_{2})x_{3}^{4}(x_{4}^{12},x_{5}^{12}) + (x_{1},x_{2})(x_{4}^{4},x_{5}^{4})(x_{4}^{12},x_{5}^{12}) +
 x_{3}(x_{4}^{4},x_{5}^{4})(x_{4}^{12},x_{5}^{12}) + \]\[ + x_{3}^{4}(x_4,x_5)(x_{4}^{12},x_{5}^{12}) +
 x_{3}^{5}(x_{4}^{12},x_{5}^{12}) + (x_{4},x_{5})(x_{4}^{4},x_{5}^{4})(x_{4}^{12},x_{5}^{12}).\]
By Proposition $2.2$, we get $I=I^{(1)}+I^{(2)}+I^{(3)}$.
\end{exm}

\begin{obs}
For any $1\leq q\leq r$ and any nonnegative integers $\gamma_{1},\ldots,\gamma_{q} \leq_{\de} \alpha_{q}$ such that
\linebreak
$\gamma_{1}+\cdots +\gamma_{i}<\alpha_{i}$, $\gamma_{1}+\cdots +\gamma_{i}<_{\de}\alpha_{q}$ for $1\leq i<q$ and $\gamma_{1}+\cdots +\gamma_{q}=\alpha_{q}$ we denote
\[I^{(q)}_{\gamma_{1},\ldots,\gamma_{q}} = \prod_{e=1}^{q}\prod_{t=0}^{s} (\mathbf{n}_{e}^{[d_{t}]})^{\gamma_{et}}.
Proposition\; 2.2\; implies:\;
I = \sum_{q=1}^{r} \sum_{ \gamma_{1},\ldots,\gamma_{q}}I^{(q)}_{\gamma_{1},\ldots,\gamma_{q}}.\]
Let $\me=(x_{1},\ldots,x_{n})\subset S$ be the irrelevant ideal of $S$. We have:
\[ (I:_{S}\me) = \bigcap_{j=1}^{n}(I:x_{j}) = \bigcap_{j=1}^{n}( (\sum_{q=1}^{r} \sum_{ \gamma_{1},\ldots,\gamma_{q}}I^{(q)}_{\gamma_{1},\ldots,\gamma_{q}}) : x_{j}) =
   \bigcap_{j=1}^{n}( \sum_{q=1}^{r} \sum_{ \gamma_{1},\ldots,\gamma_{q}} (I^{(q)}_{\gamma_{1},\ldots,\gamma_{q}} : x_{j})). \]
On the other hand, if $x_{j}\in \mathbf{n}_{p}$ for some $1\leq p\leq q$ then 
\[ J^{(q),j}_{\gamma_{1},\ldots,\gamma_{q}} := (I^{(q)}_{\gamma_{1},\ldots,\gamma_{q}} : x_{j}) = \prod_{e\neq p}^{q}\prod_{t=0}^{s} (\mathbf{n}_{e}^{[d_{t}]})^{\gamma_{et}}\mathbf{n_{p,\hat{j}}}^{[d_t]} (\mathbf{n}^{[d_t]})^{\gamma_{pt}-1}(\sum_{\gamma_{pt}>0}\prod_{j\neq t}(\mathbf{n}_{e}^{[d_{t}]})^{\gamma_{jt}}), \]
where $\mathbf{n_{p,\hat{j}}}^{[d_t]}= (x_{i_{p-1}+1}^{d_t},\ldots,x_{j}^{d_t-1},\ldots,x_{i_{p}}^{d_t})$
and $\mathbf{n_{p,\hat{j}}}^{[d_t]} (\mathbf{n}^{[d_t]})^{\gamma_{pt}-1}:=S$ if $\gamma_{pt}=0$. Thus
\[ (I:_{S}\me) =  \sum_{q^{1}=1}^{r} \sum_{ \gamma_{1}^{1},\ldots,\gamma_{q^1}^{1}} 
\cdots  \sum_{q^{n}=1}^{r} \sum_{ \gamma_{1}^{n},\ldots,\gamma_{q^n}^{n}}  \bigcap_{j=1}^{n} J^{(q^{j}),j}_{\gamma_{1}^{j},\ldots,\gamma_{q^j}^{j}}, \]
where for a given $q=q^{j}$, we take the second $j^{th}$ sum for $\gamma_{1}^{j},\ldots,\gamma_{q}^{j}\leq_{\de} \alpha_{q}$ such that $\gamma_{1}^{j}+\cdots +\gamma_{i}^{j}<\alpha_{i}$, $\gamma_{1}^{j}+\cdots +\gamma_{i}^{j}<_{\de}\alpha_{q}$ for $1\leq i<q^{j}$ and $\gamma_{1}^{j}+\cdots +\gamma_{q}^{j}=\alpha_{q}$.
\end{obs}

\begin{prop}
Let $n\geq 2$ and let $1\leq i_{1} < i_{2} < \cdots < i_{r}=n$ be some integers. Let $\alpha_{1}< \alpha_{2}< \cdots < \alpha_{r}$ be some positive integers. 
We consider the ideal $I=\sum_{q=1}^{r}I_{q}$, where $I_q = <x_{i_{q}}^{\alpha_{q}}>_{\de}$. Then, we have: $reg(I)\leq reg(I_r)$ (We will see later in which conditions we have equality).
\end{prop}

\begin{proof}
From \cite[Corollary 3.6]{mir} it follows that $(I_{q})_{\geq e}$ is stable, if $e\geq reg(I_q)$ so $(I_{q})_{\geq e}$ is stable for $e=max\{reg(I_1),\ldots,reg(I_r)\}$. Since $I_{\geq e} = \sum_{q=1}^{r}(I_{q})_{\geq e}$ and since a sum of stable ideals is still a stable ideal, it follows that $I_{\geq e}$ is stable. Therefore, from \cite[Proposition 12]{ert} we get $reg(I)\leq e$. On the other hand, if we denote $s_q=max\{t|\; \alpha_{qt}>0\}$ for any $1\leq q\leq r$, from \cite[Theorem 3.1]{mir} we get $reg(I_q) = \alpha_{qs_q}d_{s_q} + (i_q - 1)(d_{s_q}-1)$, thus $max\{reg(I_1),\ldots,reg(I_r)\}=reg(I_r)$. In conclusion, $reg(I)\leq reg(I_r)$. 
\end{proof}

\begin{prop}
With the above notations, for any $1\leq q\leq r$ we have:
\[ (I_q:\me_{q}) + (I_1 + \cdots +I_q) \subset ((I_1+\cdots+I_q):\me_{q}) \subset
((I_1+\cdots+I_q):\en_{q}) = (I_q:\en_{q}) + (I_1 + \cdots +I_q).\]
\end{prop}

\begin{proof}
Fix $1\leq q\leq r$. The first two inclusions are obvious. In order to prove the last equality, it is enough to
show that $((I_1+\cdots+I_q):x_j) \subset (I_q:x_j) + (I_1 + \cdots +I_q)$ for any $x_j\in \en_q$. Indeed, suppose 
$u\in ((I_1+\cdots+I_q):x_j)$, therefore $x_j\cdot u \in I_1+\cdots+I_q$. If $x_j\cdot u \notin I_q$ it follows that
$x_j\cdot u \in I_e$ for some $e<q$. Thus $u\in I_e$, since $x_j$ does not divide any minimal generator
of $I_e$.
\end{proof}

Let $n\geq 2$ and let $1\leq i_{1} < i_{2} < \cdots < i_{r}=n$ be some integers. Let $\alpha_{1}< \alpha_{2}< \cdots < \alpha_{r}$ be some positive integers. We write $\alpha_q=\sum_{t\geq 0}\alpha_{qt}d_t$.
Let $s_q=max\{t|\; \alpha_{qt}>0\}$ for any $1\leq q\leq r$. Notice that $s_1\leq s_2 \leq \cdots \leq s_r$. Let $1\leq q_1 < q_2 <\cdots <q_k = r$ such that:
\[ s_1=\cdots = s_{q_1} < s_{q_1 + 1} = \cdots = s_{q_2} < \cdots < s_{q_{k-1}+1} = \cdots = s_{q_k}. \]
For $1\leq j \leq k$ we define some positive integers $\chi_j$ as follows. If
$i_{q_{j}}-i_{q_{j}-1}\geq 2$ we put $\chi_{j} = (d_{s_{q_{j}}} - 1)(i_{q_{j}}-i_{q_{j-1}}) + d_{s_{q_j}}(\alpha_{q_j s_{q_j}}-1)$. Otherwise, suppose that $q=q_{j}$ and there exists a positive integer $1\leq l\leq r-q+1$ such that $s_{q-1}<s_{q}<\cdots<s_{q+l-1}$ and $i_{q+l-1}=i_{q-1}+l$. Denote $i=i_{q}$. We define recursively the numbers
$\chi_{i+m-1}$, for $1\leq m\leq l$, starting with $m=l$. Suppose that we already define $\chi_{i+m},\ldots,\chi_{i+l-1}$. If $\alpha_{q+m-2,s_{q+m-2}}>\alpha_{q+m-1,s_{q+m-1}}$, we put
$\chi_{q+m-1}:=\sum\limits_{t=s_{q+m-2}+1}^{s_{q+m-1}} \alpha_{q+m-1,t}d_t - 1$ and we switch from $m$ to $m-1$.
Otherwise, if $\alpha_{q+m-2,s_{q+m-2}}\leq \alpha_{q+m-1,s_{q+m-1}}$ we put
\[ \chi_{q+m-1}:= (\alpha_{q+m-1,s_{q+m-2}}-\alpha_{q+m-2,s_{q+m-2}}+1)\cdot d_{s_{q+m-2}} + \sum\limits_{t=s_{q+m-2}+1}^{s_{q+m-1}} \alpha_{q+m-1,t}d_t - 1 \]
and, if $m\geq 2$, we put also $\chi_{q+m-2}:=\alpha_{q+m-2,s_{q+m-2}}\cdot d_{s_{q+m-2}} - 1$. We switch from $m$
to $m-2$. We continue this procedure until $m\leq 0$. 

With these notations, for the ideal $I=<x_{i_{1}}^{\alpha_{1}},x_{i_{2}}^{\alpha_{2}},\ldots,x_{i_{r}}^{\alpha_{r}}>_{\de}$, we have the following theorem: \pagebreak

\begin{teor}
$max\{e:\;(Soc(S/I))_{e}\neq 0\} = \sum_{j=1}^{k} \chi_{j}$.
\end{teor}

\begin{proof}
For each integer $1\leq j \leq k$, we consider the following ideal:
\[ J_{j} = \begin{cases} (x_{i_{q_j}}^{\chi_j}),\;\; if\; i_{q_{j}}-i_{q_{j}-1}=1,\\
(x_{i_{q_{j-1}}+1}\cdots x_{i_{q_{j}}})^{d_{s_{q_j}}-1}\cdot \sum_{e=q_{j-1}+1}^{q_{j}} (\en_{e}^{[d_{s_{q_j}}]})^{\alpha_{es_e}-1},\;otherwise.
           \end{cases}\]
Let $J=J_1\cdot J_2 \cdots J_k$. We claim the following:\vspace{3pt}

	(1) $J\subset (I:\me)$, (2) $G(J)\cap G(I) = \emptyset$ and \vspace{3pt}
	
	(3) $max\{e|\; (Soc(S/I))_{e}\neq 0\} = max\{e|\; ((J+I)/I)_e \neq 0\}$.\vspace{3pt}

Suppose that we proved $(1),(2)$ and $(3)$. $(1)$ and $(2)$ implies $max\{e|\; ((J+I)/I)_e \neq 0\} =deg(J):=max\{deg(u)|\;u \in G(J)\}$. On the other hand, it is obvious that $deg(J)=\sum_{j=1}^{k}\chi_{j}$ and thus,
by $(3)$, we complete the proof of the theorem.

In order to prove $(1)$, we pick $x_i\in \en_{q}$ a variable, where $q\in\{1,\ldots,r\}$. Let $j$ is the unique integer with the property that $q\in \{q_{j-1}+1,\ldots,q_{j}\}$. We want to show that $x_i\cdot J\subset I$. We consider two cases.
First, we assume $i_{q_{j}}-i_{q_{j-1}} \geq 2$. We claim that $x_i J_{j} \subset I_{q_{j-1}+1} + \cdots + I_{q_{j}}$. Indeed, for any $e\in \{q_{j-1}+1,\ldots,q_{j}\}$, 
$ x_{i}(x_{i_{q_{j-1}}+1}\cdots x_{i_{q_{j}}})^{d_{s_{q_j}}-1} (\en_{e}^{[d_{s_{q_j}}]})^{\alpha_{es_e}-1} 
\subset I_{e}$,  thus $x_i J_{j} \subset I_{q_{j-1}+1} + \cdots + I_{q_{j}}$, as required. (See the proof of \cite[Lema 2.1]{mir} for details.)

Suppose now $i_{q_{j}}-i_{q_{j}-1}=1$. Let $j'\leq j$, such that if we denote $q=q_{j'}$, there exists an positive integer $j-j'+1 \leq l$ with $s_{q-1}<s_{q}<\cdots<s_{q+l-1}$, $i_{q+l-1}=i_{q-1}+l$ and $i_{q_{j'+l}}>i_{q+l-1}+1$ when
$q+l-1<r$. We prove in fact that $x_i \cdot J_{j'} \cdots J_{j} \subset I_{j}$. Note that $i=i_{q+m-1}$, where $m=j-j'+1$. Assume $m\geq 2$.
If $\alpha_{q+m-2,s_{q+m-2}} > \alpha_{q+m-1,s_{q+m-2}}$, then 
\[ x_i \cdot J_{q+m-2}J_{q+m-1} = (x_{i-1}^{\cdots + \alpha_{q+m-2,d_{s_{q+m-2}}} - 1} \cdot x_i^{\sum_{t=s_{q+m-2}+1}^{s_{q+m-1}} \alpha_{q+m-1,t}d_t}) \subset I_{j}, \]
because $\alpha_{q+m-2,d_{s_{q+m-2}}} - 1 \geq \alpha_{q+m-1,d_{s_{q+m-2}}} + d_{s_{q+m-2}} - 1$
and therefore
\[ x_i \cdot J_{q+m-2}J_{q+m-1} \subset (x_{i-1}^{d_{s_{q+m-2}} - 1} \cdot x_{i-1}^{\alpha_{q+m-1,d_{s_{q+m-2}}}} \cdot
 x_i^{\sum_{t=s_{q+m-2}+1}^{s_{q+m-1}} \alpha_{q+m-1,t}d_t}).\]
Now, the above assertion it is obvious. If $m=1$ the same trick works, with the only difference that the first "$=$" is replaced by "$\subseteq$".

If $m\geq 2$ and $\alpha_{q+m-2,s_{q+m-2}}\leq \alpha_{q+m-1,s_{q+m-2}}$ then $x_i \cdot J_{q+m-2}J_{q+m-1}$ is the ideal
\[(x_{i-1}^{\alpha_{q+m-2,d_{s_{q+m-2}}}d_{s_{q+m-2}} - 1} \cdot 
x_i^{ (\alpha_{q+m-1,s_{q+m-2}}-\alpha_{q+m-2,s_{q+m-2}}+1)d_{s_{q+m-2}} + \sum\limits_{t=s_{q+m-2}+1}^{s_{q+m-1}}
\alpha_{q+m-1,t}d_t})\]
By regrouping, we see that $x_i \cdot J_{q+m-2}J_{q+m-1} = ( x_{i-1}^{d_{s_{q+m-2}-1}} \cdot (x_{i-1}^{(\alpha_{q+m-2,d_{s_{q+m-2}}}-1)d_{s_{q+m-2}}} \cdot $

$x_i^{ (\alpha_{q+m-1,s_{q+m-2}}-\alpha_{q+m-2,s_{q+m-2}}+1)d_{s_{q+m-2}}}) \cdot x_{i}^{\sum\limits_{t=s_{q+m-2}+1}^{s_{q+m-1}}
\alpha_{q+m-1,t}d_t}) \subset I_{j}$, as required. If $m=1$ the same trick works, with the only difference that the first "$=$" is replaced by "$\subseteq$".
 
In order to prove $(2)$ it is enough to show for any $1\leq j\leq k$ that $G(J_1\cdots J_{j})\cap G(I_e) = \emptyset$
for any $e\in \{q_{j-1}+1,\ldots,q_{j}\}$, because all of the minimal generators of $J_1\cdots J_{j}$ does not contain
variables $x_i$ with $i>i_{q_j}$. We use induction on $1\leq j\leq k$. If $j=1$, then $G(J_1)\cap G(I_1) = \emptyset$
from \cite[Lemma 2.1]{mir}. Suppose the assertion is true for $j-1$. We must consider two cases.

First, suppose $i_{q_j} - i_{q_{j-1}} \geq 2$. It follows $J_j = (x_{i_{q_{j-1}}+1}\cdots x_{i_{q_{j}}})^{d_{s_{q_j}}-1}\cdot \sum_{e=q_{j-1}+1}^{q_{j}} (\en_{e}^{[d_{s_{q_j}}]})^{\alpha_{es_e}-1}$.
Since $s_{q_{j-1}}<s_{q_{j}}$ it follows that 
$J_1\cdots J_{j-1}\cdot J_j \subset (x_1,\ldots,x_{i_{q_{j-1}}})^{d_{s_{q_j}}-1}J_j$,  
and it is easy to note that none of the minimal generator of the ideal from left is included in some $I_{e}$ with
$q_{j-1}+1 \leq e \leq q_{j}$. 

Suppose now $i_{q_j} - i_{q_{j-1}} = 1$. Let $j'\leq j$, such that if we denote $q=q_{j'}$, there exists an positive integer $j-j'+1 \leq l$ with $s_{q-1}<s_{q}<\cdots<s_{q+l-1}$, $i_{q+l-1}=i_{q-1}+l$ and $i_{q_{j'+l}}>i_{q+l-1}+1$ when
$q+l-1<r$. We prove in fact that $x_i \cdot J_{j'} \cdots J_{j} \subset I_{j}$. Note that $i=i_{q+m-1}$, where $m=j-j'+1$. Assume $m\geq 2$. If $\alpha_{q+m-2,s_{q+m-2}} > \alpha_{q+m-1,s_{q+m-2}}$, then 
\[ J_1\cdots J_{j} = (J_1\cdots J_{j-2})\cdot (x_{i-1}^{\cdots + \alpha_{q+m-2,d_{s_{q+m-2}}} - 1} \cdot x_i^{\sum_{t=s_{q+m-2}+1}^{s_{q+m-1}} \alpha_{q+m-1,t}d_t -1}) \subset \]
\[  (x_1,\ldots,x_{i_{q_{j-2}}})^{d_{s_{q_{j-1}}}-1} (x_{i-1}^{\cdots + \alpha_{q+m-2,d_{s_{q+m-2}}} - 1} \cdot x_i^{\sum_{t=s_{q+m-2}+1}^{s_{q+m-1}} \alpha_{q+m-1,t}d_t -1}), \]
and it is easy to see that none of the minimal generators of the last ideals are in $I_{j}$. The subcase
$\alpha_{q+m-2,s_{q+m-2}} \leq \alpha_{q+m-1,s_{q+m-2}}$ is similar. Also, the case $m=1$.

In order to prove $(3)$ it is enough to show the "$\leq$" inequality, since obviously \linebreak $(J+I)/I\subset Soc(S/I)$. Let $u = x_{1}^{\beta_1}\cdots x_{n}^{\beta_n} \in (I:\me)$ be a monomial such that $u\notin I$. We claim that
$deg(u)\leq \sum\limits_{j=0}^{k}\chi_j$. More precisely, we claim the following:

\noindent
(a) $\sum_{i=i_{q_{j-1}}+1}^{i_{q_{j}}}\beta_i \leq \chi_j$, for all $1\leq j\leq r$ such that $i_{q_{j}}-i_{q_{j-1}}\geq 2$,

\noindent
(b) For each $j$ with the property that there exists an positive integer $1\leq l\leq r-q+1$ (where $q=q_j$) such that $s_{q-1}<s_{q}<\cdots<s_{q+l-1}$, $i_{q_{j}}-i_{q_{j-1}}\geq 2$ and $i_{q+l-1}=i_{q-1}+l$, we have $\sum_{i=i_{q_{j-1}+1}}^{i_{q_{j-1+l}}} \beta_i\leq \sum_{m=1}^{l}\chi_{j+m-1}$.

Obviously, $(a)$ and $(b)$ implies $(3)$. In order to prove (a), assume that 
$\sum_{i=i_{q_{j-1}}+1}^{i_{q_{j}}}\beta_i > \chi_j$, therefore $\sum_{i=i_{q_{j-1}}+1}^{i_{q_{j}}}\beta_i \geq 
(d_{s_{q_{j}}} - 1)(i_{q_{j}}-i_{q_{j-1}}-1) + \alpha_{q_j s_{q_j}}d_{s_{q_j}}$. It follows that we can write
$u_j = x_{i}^{d_{s_{q_j}}-1}\cdot w$, with 
$w\in (x_{i_{q_{j-1}}+1}^{d_{s_{q_j}}},\ldots,x_{i_{q_{j}}}^{d_{s_{q_j}}})^{\alpha_{q_j s_{q_j}}}$, for some
$i\in \{x_{i_{q_{j-1}}+1},\ldots, x_{i_{q_{j}}} \}$, and thus $u_j\in I_{q_{j}}$, a contradiction.
Consider now the case (b) and assume that 
$\sum_{i=i_{q_{j-1}+1}}^{i_{q_{j-1+l}}} \beta_i > \sum_{m=1}^{l}\chi_{j+m-1}$. Using similar arguments as in the case (a), we get $u_j\in I_{q_{j}}$, a contradiction.
\end{proof}

\begin{cor}
With the previous notations, $reg(I)=\sum_{j=1}^{k} \chi_{k}+1$.
\end{cor}

\begin{proof}
Since $I$ is an artinian ideal, $reg(I)=max\{e:\;Soc(S/I)_{e}\neq 0\}+1$ so the required result follows immediately from the previous theorem.
\end{proof}

\begin{obs}
We already seen that $reg(I)\leq reg(I_r)$. Now, we are able to say when we have equality, and this is only in the
case when $k=1$, i.e. $s_1 = s_2 =\cdots = s_r$. Indeed, if $k=1$, by \cite[3.1]{mir}, $reg(I_r) = (d_{s_{r}} - 1)(n-1) + d_{s_{r}}(\alpha_{r s_{r}}-1) + 1 = \chi_1 + 1$. Conversely, if $k>1$ then $\chi_1+\cdots+\chi_k < reg(I_r)$,
because $\chi_j < (d_{s_{r}} - 1)(i_{q_{j}}-i_{q_{j}-1}) + d_{s_{r}}(\alpha_{r s_{r}}-1)$ for any $j<k$.
\end{obs}

\begin{exm}
\begin{enumerate}
	\item Let $\de:1|2|6|12$ and $I=<x_{2}^{7},x_{3}^{10},x_{5}^{17}>_{\de} \subset K[x_{1},\ldots,x_{5}]$. 
	      We have $k=2$, $\chi_1 = 15$ and $\chi_2 = 22$. Therefore, $reg(I)=27$. An element of maximal degree
	      in $Soc(S/I)$ is $x_{1}^{5}x_{2}^{5}x_{3}^{5}x_{4}^{11}x_{5}^{11}$.	      
  \item Let $\de:1|4|12$ and $I=<x_1^{2},x_2^7,x_3^{16}>_{\de} \subset K[x_1,x_2,x_3]$. We have $k=3$. Since
	      $2=2\cdot 1$, $7=3\cdot 1 + 1\cdot 4$ and $16 = 1\cdot 4 + 1\cdot 12$, we get
	      $\chi_1 = 1$, $\chi_2 = 3$ and $\chi_3=19$. Therefore, $reg(I)=23$. An element of maximal degree
	      in $Soc(S/I)$ is $x_{1}x_{2}^{3}x_{3}^{19}$.
\end{enumerate}
\end{exm}

\vspace{2mm} \noindent {\footnotesize
\begin{minipage}[b]{10cm}
 Mircea Cimpoea\c s, Junior Researcher\\
 Institute of Mathematics of the Romanian Academy\\
 Bucharest, Romania\\
 E-mail: mircea.cimpoeas@imar.ro
\end{minipage}}

\begin{thebibliography}{99}
  \bibitem[1]{ah}Annetta Aramova, J\"urgen Herzog "p-Borel principal ideals", Illinois J.Math.41,no 1.(1997),103-121.
  \bibitem[2]{bs} D. Bayer, M. Stillman "A criterion for detecting m-regularity", Invent. Math 87 (1987) 1-11.
  \bibitem[3]{cs} G. Caviglia, E. Sbarra, "Characteristic-free bounds for the Castelnuovo Mumford regularity", Compos.
              Math. 141(2005), no.6, 1365-1373.
  \bibitem[4]{mir}Mircea Cimpoea\c s "A generalisation of Pardue's formula", 
             Bull. Math. Soc. Sci. Math. Roumanie (N.S.) 49(97), no. 4, 2006.             
  \bibitem[5]{mir2}Mircea Cimpoea\c s "A stable property of Borel type ideals", to appear in Communications in Algebra.    \bibitem[6]{ert}D.Eisenbud, A.Reeves, B.Totaro "Initial ideals, veronese subrings and rates of algebras", Adv.Math.   
             109 (1994), 168-187.
  \bibitem[7]{hpv}J\"urgen Herzog, Dorin Popescu, Marius Vladoiu "On the Ext-Modules of ideals of Borel type",
             Contemporary Math. 331 (2003), 171-186.
  \bibitem[8]{hp}J\"urgen Herzog, Dorin Popescu "On the regularity of p-Borel ideals", Proceed.of AMS, Volume 129,
             no.9, 2563-2570.  
  \bibitem[9]{saf}Anwar Imran, Ahmad Sarfraz "Regularity of ideals of Borel type is linearly bounded", Preprint, 2006.
\end{thebibliography}
\end{document}